\documentclass{article}

\usepackage{amsmath}
\usepackage{amssymb}
\usepackage{latexsym}
\usepackage{theorem}
\usepackage[all,2cell]{xypic}

\usepackage{graphicx}

{\theorembodyfont{\rmfamily}
\newtheorem{Definition}{Definition}[section]
\newtheorem{Example}{Example}[section]

\newtheorem{Remark}{Remark}[section]
}
\newtheorem{Theorem}{Theorem}[section]
\newtheorem{Proposition}[Theorem]{Proposition}
\newtheorem{Lemma}[Theorem]{Lemma}
\newtheorem{Corollary}[Theorem]{Corollary}
\newenvironment{Proof}{\noindent{\bf Proof\ }}{\mbox{}\hfill$\Box$\medskip}

\newcommand{\Id}{\mathrm{Id}}
\newcommand{\Frm}[1]{\mathbf{Frm}_{#1}}
\newcommand{\DFrm}[1]{\mathbf{DFrm}_{#1}}
\newcommand{\PD}[1]{\mathbf{PD}_{#1}}

\newcommand{\dom}{\mathrm{dom}}
\newcommand{\cod}{\mathrm{cod}}

\begin{document}

\UseAllTwocells

\title{Weak $\omega$-Categories as $\omega$-Hypergraphs}
\author{Hiroyuki Miyoshi \\
{\normalsize Department of Computer Science, Kyoto Sangyo University} \\
{\normalsize\tt hxm@cc.kyoto-su.ac.jp}
\and
Toru Tsujishita\\
{\normalsize Department of Mathematics, Hokkaido University}\\
{\normalsize\tt tujisita@math.sci.hokudai.ac.jp}
}
\maketitle

\begin{abstract}
In this paper, firstly, we introduce a higher-dimensional
analogue of hypergraphs, namely {\em $\omega$-hypergraphs}.
This notion is thoroughly flexible because
unlike ordinary $\omega$-graphs,
an $n$-dimensional edge called an {\em $n$-cell}
has many sources and targets. Moreover, cells have polarity,
with which pasting of cells is implicitly defined.
As examples, we also give some known structures in terms of
$\omega$-hypergraphs. 
Then we specify a special type of $\omega$-hypergraph,
namely {\em directed $\omega$-hypergraphs},
which are made of cells with direction.
Finally, besed on them, we construct our {\em weak $\omega$-categories}.
It is an $\omega$-dimensional variant of
the weak $n$-categoreis given by Baez and Dolan \cite{BD}.
We introduce {\em $\omega$-identical},
{\em $\omega$-invertible} and  {\em $\omega$-universal} cells
instead of universality and balancedness in \cite{BD}.
The whole process of our definition is in parallel with the way of
regarding categories as graphs with composition and identities.
\end{abstract}

\tableofcontents

\section{Introduction}

J. Baez and J. Dolan recently proposed an important
and impressive definition of weak $n$-categories\cite{BD}.
They utilize nonstandard $n$-cells with not just one but many
$n-1$-cells as their domains for taming coherence conditions.
Authors' primary motivation was to understand their idea
along the famous slogan ``categories are graphs with monoid
structures''.
Thus they investigated a suitable notion of $n$- or
$\omega$-dimensional graph-like structures
which should include the underlying structures
of Baez-Dolan-style weak $\omega$-categories.

In the way of pursuing such structures,
they found a general notion of $\omega$-dimensional structures
whose $n$-cells have many $n-1$-cells not only in their domains
but also in their codomains. This notion contains
various categorical algebras:
$\omega$-categories, bicategoreis, double categories, etc.
Meanwhile, authors noticed that it can be thought of
as a form of $\omega$-dimensional hypergraphs.
Hypergraphs have been explored in mathematics\cite{CB},
database theory\cite{RF}, concurrency theory\cite{DH} and
graph rewriting\cite{DP} as a device to represent complex notions.
But their higher-dimensional extensions are still not known
corresponding to $n$- or $\omega$-graphs for ordinary graphs.
Therefore such structures are named {\em $\omega$-hypergraphs}\footnote{The  
definition of $\omega$-hypergraphs in this paper
is not the most general form,
because each node is shared by at most two hyperedges.}.

Thus the purpose of this paper is two-folded:
One is to provide a general environment
for representing various concepts,
especially developing various category theories.
Another is to give a definition of weak $\omega$-categories
which respects saturatedness in the meaning of M. Makkai\cite{MM1}.

\section{Trees and forests}

\begin{figure}
  \begin{center}
  \includegraphics{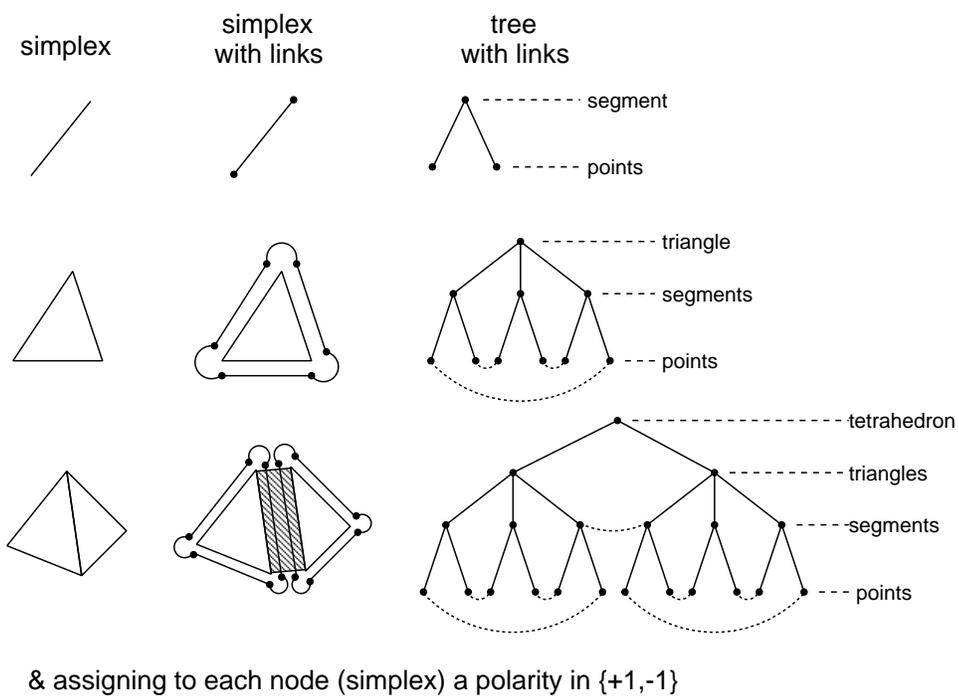}
  \end{center}
  \caption{From simplices to trees with links and polarity\label{motivation:fig1}}
 \end{figure}

Our main idea is to represent an $n$-cell as a tree with links and polarity.
This is refinement of usual simplice (Figure \ref{motivation:fig1}).
We start with the definition of trees and forests.

\begin{Definition}[$n$-trees and $n$-forests]
For any natural number $n \geq 0$, an {\em $n$-tree} $T$ is a triple
$\langle r^T,\ S^T,\ \pi^T \rangle$ consisting of
\begin{itemize}

\item $S^T_n = \{r^T\}$, whose element $r^T$ is called the {\em root} of $T$;

\item $S^T = \coprod_{0\leq i \leq n} S^T_i$,
where $S^T_i$ is a finite set whose elements are called
{\em $i$-nodes} or simply {\em nodes};

\item $\pi^T =\coprod_{0\leq i \leq n-1} \pi^T_i$,
where $\pi^T_i$ is a function from $S^T_i$ to $S^T_{i+1}$.

\end{itemize}
Also, an {\em $n$-forest} $F$ is a pair $\langle S^F,\ \pi^F \rangle$,
consisting of
\begin{itemize}

\item $S^F = \coprod_{0\leq i \leq n} S^F_i$,
where $S^F_i$ is a finite set of {\em $i$-nodes};

\item $\pi^F =\coprod_{0\leq i \leq n-1} \pi^F_i$,
where $\pi^F_i$ is a function from $S^F_i$ to $S^F_{i+1}$.

\end{itemize}
\end{Definition}

\begin{Definition}[isomorphism of trees and forests]
For any natural number $n\geq 0$,
a {\em homomorphism} of $n$-trees $\sigma: T \rightarrow T'$
is a map from $S^T$ to $S^{T'}$ such that, for every $x \in S^T_i$,
$\sigma(x)\in S^{T'}_i$ and 
$\sigma \circ \pi^T = \pi^{T'}\circ \sigma$.
A homomorphism $\sigma$ is an {\em isomorphism} when it is a bijection.
A {\em homomorphism} and an {\em isomorphism} of $n$-forests
are also defined in the same way.
\end{Definition}

\begin{figure}
 \begin{center}
  \includegraphics{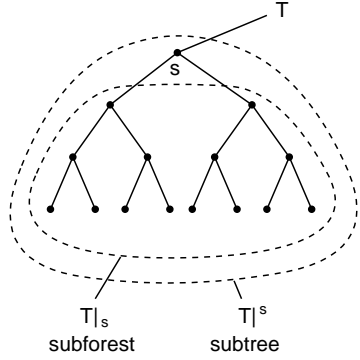}
 \end{center}
 \caption{the subtree and subforest at $s$\label{subtree:fig1}}
\end{figure}
\begin{Definition}[subtrees and subforests]
For an $n$-tree $T=\langle r^T,\ S^T,\ \pi^T \rangle$ 
and a $k$-node $s$ ($0\leq k\leq n$),
a {\em subtree} with the root $s$ is defined as
$T|^s = \langle s,\ S^{T|^s}, \pi^{T|^s} \rangle$ where
\begin{itemize}

\item $S^{T|^s} = \coprod_{0\leq i \leq k} S^{T|^s}_i$,
	where $S^{T|^s}_i = \{t\in S^{T}_i\,|\,(\pi^{T})^{k-i}(t)=s\}$;

\item $\pi^{T|^s} = \coprod_{0\leq i \leq k-1} \pi^{T|^s}_i$
	where $\pi^{T|^s}_i = \pi_i|_{S^{T|^s}_i}$

\end{itemize}
And for a $k$-node $s$ ($1\leq k\leq n$),
a {\em subforest} under $s$ is defined as
$T|_s = \langle S^{T|_s}, \pi^{T|_s} \rangle$ where
\begin{itemize}

\item $S^{T|_s} = \coprod_{0\leq i \leq k-1} S^{T|_s}_i$,
	where $S^{T|_s}_i = \{t\in S^{T}_i\,|\,(\pi^{T})^{k-i}(t)=s\}$;

\item $\pi^{T|_s} = \coprod_{0\leq i \leq k-2} \pi^{T|_s}_i$
	where $\pi^{T|_s}_i = \pi_i|_{S^{T|_s}_i}$

\end{itemize}

Also in the same way,
for an $n$-forest $F=\langle S^F,\ \pi^F \rangle$
and a $k$-node $s$ ($0\leq k\leq n$),
a {\em subtree} with the root $s$ is defined as
$F|^s = \langle s,\ S^{F|^s}, \pi^{F|^s} \rangle$,
and for a $k$-node $s$ ($1\leq k\leq n$),
a {\em subforest} under $s$ as
$F|_s = \langle S^{F|_s}, \pi^{F|_s} \rangle$.
\end{Definition}

\section{Shells}

\begin{figure}
 \begin{center}
  \includegraphics{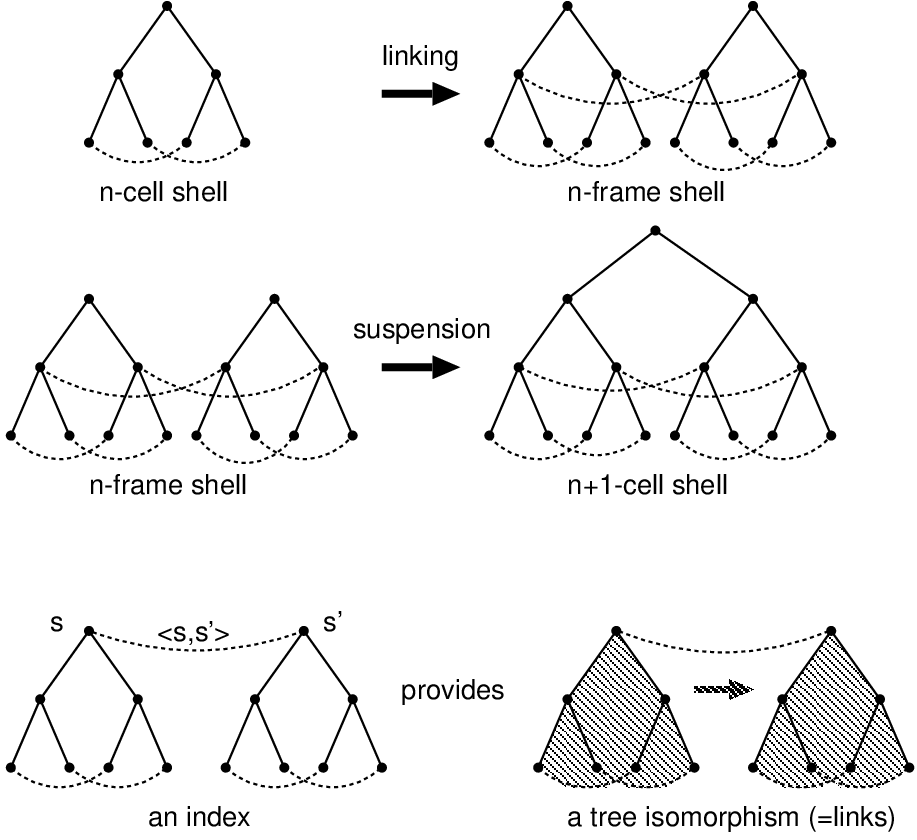}
 \end{center}
 \caption{Mutually inductive construction of shells\label{shell:fig1}}
\end{figure}

Shells play the same role as shape diagrams
in the ordinary category theory.
We will mutually inductively define a shell for each cell
as a tree with polarity and links
and one for each frame as a forest with polarity and links
(Figure \ref{shell:fig1}).

\subsection{the base case}

\begin{Definition}[$0$-cell shells and $0$-frame shells]
A {\em $0$-cell shell} $\theta$ is a singleton set with polarity.
More precisely, it is
$\langle r^{\theta}, S^{\theta}, \emptyset, \epsilon^{\theta}, \emptyset, \emptyset \rangle$
where $S^{\theta}=S^{\theta}_0=\{r^{\theta}\}$
and $\epsilon^{\theta}$ is a function from $S^{\theta}$ to $\{-1,1\}$. 
Similarly, A {\em $0$-frame shell} $\xi$ is a set with polarity, that is, 
$\langle S^{\xi}, \emptyset, \epsilon^{\xi}, \emptyset, \emptyset \rangle$ where $S^{\xi}=S^{\xi}_0$
and $\epsilon^{\xi}$ is a function from $S^{\xi}$ to $\{-1,1\}$. 
\end{Definition}
For the meaning of these definitions, see the following sections.

\subsection{the induction step}

Suppose that for the dimensions less than $n$,
all staff has already been defined.

\begin{Definition}[$n$-cell shells]
An {\em $n$-cell shell} $\theta$ is
\[
\langle r^{\theta},
\ S^{\theta},
\ \pi^{\theta},
\ \epsilon^{\theta},
\ \Upsilon^{\theta},
\ \{\sigma^{\theta}_{\langle s,s' \rangle}\}_{{\langle s,s' \rangle} \in \Upsilon^{\theta}}
\rangle
\]
consisting of the following data:
\begin{itemize}

\item $\underline{\theta}= 
\langle r^{\theta},\ S^{\theta},\ \pi^{\theta} \rangle$ is an $n$-tree,
called the {\em base $n$-tree} of $\theta$;

\item {\em polarity:}
$\epsilon^{\theta}$ is a function from $S^{\theta}$ to $\{-1,1\}$;

\item {\em links:}
$\Upsilon^{\theta}=\coprod_{0\leq i \leq n-2} \Upsilon^{\theta}_i$,
where $\Upsilon^{\theta}_i \subset S^{\theta}_i \times S^{\theta}_i$;

\item {\em linking isomorphisms:}
for each $\langle s,s' \rangle \in \Upsilon^\theta_i$,
$\sigma_{\langle s,s' \rangle}$ is a $i$-tree isomorphism
from $\underline{\theta}|^s$ to $\underline{\theta}|^{s'}$,

\end{itemize}
which satisfy the following condition:
\begin{itemize}

\item {\em mutuality:} $\theta|_{r^{\theta}} = 
\langle 
  S^{\theta|_{r^{\theta}}},
\ \pi^{\theta|_{r^{\theta}}},
\ \epsilon^{\theta|_{r^{\theta}}},
\ \Upsilon^{\theta},
\ \{\sigma^{\theta}_{\langle s,s' \rangle}\}_{{\langle s,s' \rangle} \in \Upsilon^{\theta}}
\rangle$ is an $n-1$-frame shell,
where
  \begin{itemize}

  \item $\langle S^{\theta|_{r^{\theta}}},\ \pi^{\theta|_{r^{\theta}}} \rangle
		= \underline{\theta}|_{r^{\theta}}$

  \item $\epsilon^{\theta|_{r^{\theta}}} = \epsilon^{\theta}|_{S^{\theta|_{r^{\theta}}}}$

  \end{itemize}
\end{itemize}

\end{Definition}

\begin{Definition}[$n$-frame shells]
An {\em $n$-frame shell} $\xi$ is 
\[
\langle 
  S^{\xi},
\ \pi^{\xi},
\ \epsilon^{\xi},
\ \Upsilon^{\xi},
\ \{\sigma^{\xi}_{\langle s,s' \rangle}\}_{{\langle s,s' \rangle} \in \Upsilon^{\xi}}
\rangle
\]
consisting of the following data:
\begin{itemize}

\item $\underline{\xi} = \langle S^{\xi},\ \pi^{\xi} \rangle$ is an $n$-forest,
called the {\em base $n$-forest} of $\xi$;

\item {\em polarity:}
$\epsilon^{\xi}$ is a function from $S^{\xi}$ to $\{-1,1\}$;

\item {\em links:}
$\Upsilon^{\xi}=\coprod_{0\leq i \leq n-1} \Upsilon^{\xi}_i$,
where $\Upsilon^{\xi}_i \subset S^{\xi}_i \times S^{\xi}_i$;

\item {\em linking isomorphisms:}
for $\langle s,s' \rangle \in\Upsilon^{\xi}_i$,
$\sigma_{\langle s,s' \rangle}$ is a $i$-tree isomorphism
from $\underline{\xi}|^s$ to $\underline{\xi}|^{s'}$,

\end{itemize}
which satisfy the following conditions:
\begin{itemize}

\item {\em mutuality:} for any $s \in S^{\xi}_n$, $\xi|^s =
\langle 
s,
\ S^{\xi|^s},
\ \pi^{\xi|^s},
\ \epsilon^{\xi|^s},
\ \Upsilon^{\xi|^s},
\ \{\sigma^{\xi}_{\langle s,s' \rangle}\}_{{\langle s,s' \rangle} \in \Upsilon^{\xi|^s}}
\rangle$ is an $n$-cell shell,
where
  \begin{itemize}

  \item $\langle s,\ S^{\xi|^s},\ \pi^{\xi|^s}\rangle
	=\underline{\xi}|^s$

  \item $\epsilon^{\xi|^s}= \epsilon^{\xi}|_{S^{\xi|^s}}$;

  \item $\Upsilon^{\xi|^s} = \coprod_{0\leq i \leq n-2} \Upsilon^{\xi|^s}_i$,
        where $\Upsilon^{\xi|^s}_i
                = \Upsilon^{\xi}_i|_{S^{\xi|^s}_i \times S^{\xi|^s}_i}$

  \end{itemize}

\item {\em bijectivity:} 
if $\langle s,t \rangle, \langle s,t' \rangle \in \Upsilon^{\xi}_{n-1}$,
then $t=t'$; and
if $\langle s,t \rangle, \langle s',t \rangle \in \Upsilon^{\xi}_{n-1}$,
then $s=s'$;

\item {\em involution:}
if $\langle s, t \rangle \in \Upsilon^{\xi}_{n-1}$, then
$\langle t, s \rangle \in \Upsilon^{\xi}_{n-1}$
and
$\sigma^{\xi}_{\langle s, t \rangle}
\circ \sigma^{\xi}_{\langle t, s \rangle} = \Id$ and
$\sigma^{\xi}_{\langle t, s \rangle}
\circ \sigma^{\xi}_{\langle s, t \rangle} = \Id$;

\item {\em conjugation:}
if $\langle s,s'\rangle \in \Upsilon^{\xi}_{n-1}$ and $t \in S^{\xi}$
such that $(\pi^{\xi})^i(t)=s$ for some $i\geq 0$,
then $\epsilon(t)\epsilon(\sigma^{\xi}_{\langle s, s' \rangle}(t)) = -1$
(this implies
{\em anti-reflexivity}: if $\langle s, s' \rangle \in
\Upsilon^{\xi}_{n-1}$, then $s \neq s'$);

\item {\em correspondence of links:}
if $\langle s, s'\rangle \in \Upsilon^{\xi}_{n-1}$ and
$\langle t, t'\rangle \in \Upsilon^{\xi}_i$ for some $i\leq n-3$
and $(\pi^{\xi})^{n-i}(t)=(\pi^{\xi})^{n-i}(t')=s$,
then $\langle \sigma^{\xi}_{\langle s, s'\rangle}(t),
\sigma^{\xi}_{\langle s, s'\rangle}(t') \rangle \in \Upsilon^{\xi}_i$
and
$\sigma^{\xi}_{\langle t, t'\rangle}\circ
\sigma^{\xi}_{\langle s, s'\rangle}=
\sigma^{\xi}_{\langle \sigma^{\xi}_{\langle s, s'\rangle}(t),
\sigma^{\xi}_{\langle s, s'\rangle}(t') \rangle}\circ
\sigma^{\xi}_{\langle t, t'\rangle}$;

\item {\em commutativity of links:}
for $k\geq 2$ and
$\langle s_1',s_2\rangle$,
$\langle s_2',s_3\rangle$,\ldots,
$\langle s_{k-1}',s_k\rangle$,
$\langle s_k',s_1\rangle$
in $\Upsilon^{\xi}$
such that $\pi^{\xi}(s_i)=s_i'$ or $\pi^{\xi}(s_i')=s_i$,
if $s'=
 (\sigma^{\xi}_{\langle s_k',s_1 \rangle} \circ
  \sigma^{\xi}_{\langle s_{k-1}',s_k \rangle} \circ \cdots \circ
  \sigma^{\xi}_{\langle s_{2}',s_3 \rangle} \circ
  \sigma^{\xi}_{\langle s_{1}',s_2 \rangle})(s)$
is defined, then $s'=s$;
\begin{figure}
 \begin{center}
  \includegraphics{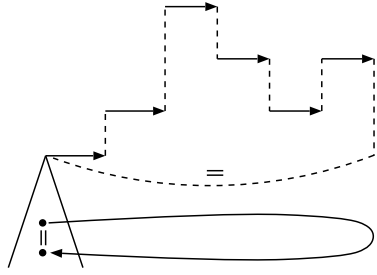}
 \end{center}
 \caption{commutativity of links\label{commlinks:fig1}}
\end{figure}
(that is,
if $s_i'$ is of the smallest level between $s_1'$,\ldots,$s_k'$, then
the composition of such isomorphisms as above
$\sigma^{\xi}_{\langle s_q',s_i \rangle} \circ \cdots \circ
\sigma^{\xi}_{\langle s_{i}',s_{p} \rangle}$ is defined
and is the identity homomorphism of the subtree at $s_{i}'$).

\item {\em closedness:}
for every $s \in S^{\xi}_{n-1}$, there exists a(n unique) node
$s' \in S^{\xi}_{n-1}$
such that $\langle s, s'\rangle \in \Upsilon^{\xi}_{n-1}$.

\end{itemize}
\end{Definition}

Closedness means {\em globularity} of higher dimensional cells.
Note that every $t \in S^{\xi}$ is
in $S^{\xi|^s}$ for just one $s \in S^{\xi}_n$;
and also every $\langle t, t'\rangle \in \Upsilon^{\xi}_i$ 
for $i\leq n-2$ is in
$\Upsilon^{\xi|^s}_i$ for just one $s \in S^{\xi}_n$.
The latter is due to the closedness of frame shells at lower levels.

\begin{Proposition}
For a cell shell $\theta$,
if $\langle s,s' \rangle \in \Upsilon^{\theta}$, then 
$\pi^i(s)=\pi^i(s')$ where $i=1$ or $2$.
\end{Proposition}

\begin{Remark}
Thus the situation of the correspondence of links for $k=2$ occur 
only when $l_1 = l_2 = 1$ and either
the parents of $s_1'$ and $s_2$ or those of $s_2$ and $s_1'$ 
are the same.
And for an $n$-cell shell $\theta$ in $n$-frame shell,
an {\em outer link} $\langle s,s' \rangle$ of which
$s$ or $s'$ is not in $S_{\theta}$, must be an $n-1$-link.
\end{Remark}

\begin{Definition}[$\cong_n$, $(-)^*$]
For two $n$-frame shells
\begin{align*}
\xi &=
\langle 
  S^{\xi},
\ \pi^{\xi},
\ \epsilon^{\xi},
\ \Upsilon^{\xi},
\ \{\sigma^{\xi}_{\langle s,s' \rangle}\}_{{\langle s,s' \rangle} \in \Upsilon^{\xi}}
\rangle\quad\text{and}  \\
\xi' &=
\langle 
  S^{\xi'},
\ \pi^{\xi'},
\ \epsilon^{\xi'},
\ \Upsilon^{\xi'},
\ \{\sigma^{\xi'}_{\langle s,s' \rangle}\}_{{\langle s,s' \rangle} \in \Upsilon^{\xi'}}
\rangle\,\text{,}
\end{align*}
an {\em isomorphism} $f$
from $\xi$ to $\xi'$ is an $n$-forest isomorphism
(with its inverse $f^{-1}$) such that
\begin{itemize}

\item $\epsilon^{\xi}(s) = \epsilon^{\xi'}(f(s))$
\quad ($\Leftrightarrow$
$\epsilon^{\xi'}(t) = \epsilon^{\xi}(f^{-1}(t))$);

\item if $\langle s,s' \rangle \in \Upsilon^{\xi}$, then
$\langle f(s),f(s') \rangle \in \Upsilon^{\xi'}$
\quad ($\Leftrightarrow$
if $\langle t,t' \rangle \in \Upsilon^{\xi'}$, then
$\langle f^{-1}(t),f^{-1}(t') \rangle \in \Upsilon^{\xi}$);

\item $f\circ\sigma^{\xi}_{\langle s,s' \rangle} =
        \sigma^{\xi'}_{\langle f(s),f(s') \rangle}\circ f$
\quad ($\Leftrightarrow$
$f^{-1}\circ\sigma^{\xi'}_{\langle t,t' \rangle} =
\sigma^{\xi}_{\langle f(t),f(t') \rangle}\circ f^{-1}$).

\end{itemize}
When an isomorphism $f$ from $\xi$ to $\xi'$ exists,
we say that $\xi$ is {\em isomorphic} to $\xi'$, and write
$f:\xi \cong_{n} \xi'$, $\xi \cong_{n} \xi'$, or simply $\xi \cong \xi'$.
Obviously $\cong_{n}$ is an equivalence relation.


For an $n$-frame shell $\xi$, $(\xi)^*$ is defined as
$\langle 
  S^{\xi},
\ \pi^{\xi},
\ \epsilon^{(\xi)^*},
\ \Upsilon^{\xi},
\ \{\sigma^{\xi}_{\langle s,s' \rangle}\}_{{\langle s,s' \rangle} \in \Upsilon^{\xi}},
\rangle$ where $\epsilon^{(\xi)^*}(s) = - \epsilon^{\xi}(s)$.
It is easy to check well-definedness,
that is, $(\xi)^*$ is in fact an $n$-frame shell, and $((\xi)^*)^*=\xi$.
\end{Definition}
An $n$-cell shell can be seen as an $n$-frame shell.
Thus we can define isomorphsims between $n$-cell shells similarly.

\section{Diagrams and $\omega$-hypergraphs}

Cell diagrams and frames are mutually inductively defined.

\begin{Definition}[$i$-cell]
We prepare a set of $i$-cells for each $i\in\mathbb{N}\cup\{0\}$:
\begin{itemize}


\item $\Sigma_i = \Sigma_{i,-1}\amalg \Sigma_{i,1}$,

\item a bijection $(-)^*:\Sigma_i \rightarrow\Sigma_i$
such that for $c\in \Sigma_{i,k}$ with $k\in\{1,-1\}$,
$c^* \in \Sigma_{i,-k}$ and $(c^*)^* = c$.

\end{itemize}
Elements of $\Sigma_{i}$ are called {\em $i$-cells};
those of $\Sigma_{i,1}$ {\em positive} $i$-cells;
and those of $\Sigma_{i,-1}$ {\em negative} $i$-cells.
$c^*$ is called the {\em conjugate} of $c$.
\end{Definition}

\subsection{the base case}

\begin{Definition}[$0$-hypergraph, $0$-cell diagram and $0$-frame]
For consistency, let $\Frm{-1}$ be $\{ \emptyset\}$,
the only one $-1$-frame isomorphism the empty function
$\emptyset:\emptyset \rightarrow \emptyset$
and $\partial_0:\Sigma_0\rightarrow\Frm{-1}$ the unique function.
A {\em $0$-hypergraph} is $\langle \Sigma_{0}, \partial_0 \rangle$.
A {\em $0$-cell diagram} $\eta$ is 
$\langle r^{\eta}, S^{\eta}, \emptyset, \epsilon^{\eta}, \emptyset, \emptyset, \lambda^{\eta},
\{\rho^{\eta}_{r^{\eta}}\}_{r^{\eta}\in S^{\eta}} \rangle$,
where $\underline{\eta}=\langle r^{\eta}, S^{\eta}, \emptyset, \epsilon^{\eta}, \emptyset, \emptyset \rangle$ is
a $0$-cell shell, $\lambda^{\eta}$ is a function from $S^{\eta}$ to $\Sigma_0$
such that $\lambda^{\eta}(r^{\eta}) \in \Sigma_{0,\epsilon^{\eta}(r^{\eta})}$, and
$\rho^{\eta}_{r^{\eta}}$ is the empty function.
A {\em $0$-frame diagram}, or simply a {\em $0$-frame}, $\zeta$ is 
$\langle S^{\zeta}, \emptyset, \epsilon^{\zeta}, \emptyset, \emptyset, \lambda^{\zeta},
\{\rho^{\zeta}_{s}\}_{s\in S^{\zeta}}\rangle$,
where $\underline{\zeta}=\langle S^{\zeta}, \emptyset, \epsilon^{\zeta}, \emptyset, \emptyset \rangle$ is
a $0$-frame shell, $\lambda^{\zeta}$ is a function from $S^{\zeta}$ to $\Sigma_0$
such that for any $s\in S^{\zeta}$, $\lambda^{\zeta}(s) \in \Sigma_{0,\epsilon^{\zeta}(s)}$, and
each $\rho^{\zeta}_{s}$ is the empty function.
A {\em $0$-frame isomorphism} from $\zeta$ to $\zeta'$ is
a $0$-frame shell isomorphism 
$f: \underline{\zeta}\rightarrow \underline{\zeta'}$
(in fact, a bijection from $S^{\zeta}$ to $S^{\zeta'}$) satisfying
$\lambda^{\zeta}=\lambda^{\zeta'}\circ f$.
$\Frm{0}$ is the set of $0$-frames.
\end{Definition}

\subsection{the induction step}

Suppose that $n\geq 1$ and that for the dimensions less than $n$,
all staff has already been defined.

\begin{Definition}[boundary of $n$-cells]
As a parameter of definitions, 
a function $\partial_{n}:\Sigma_{n}\rightarrow\Frm{n-1}$ 
satisfying $(\partial_{n}(c))^* = \partial_{n}(c^*)$
are given (for the $n-1$ dimension,
$(-)^*$ for frames have been defined).
$\partial_{n}(c)$ is called the {\em boundary} of $c$.
\end{Definition}

\begin{Definition}[$n$-hypergraph]
An {\em $n$-hypergraph} $G = \langle \Sigma, \partial \rangle$
consists of
\begin{itemize}

\item $\Sigma = \coprod_{0\leq i\leq n} \Sigma_{i}$, and

\item $\partial = \coprod_{1\leq i\leq n} \partial_{i}$.

\end{itemize}
\end{Definition}

\begin{Definition}[$n$-cell diagram]

An {\em $n$-cell diagram} $\eta$ is
\[
\langle r^{\eta},
\ S^{\eta},
\ \pi^{\eta},
\ \epsilon^{\eta},
\ \Upsilon^{\eta},
\ \{\sigma^{\eta}_{\langle s,s' \rangle}\}_{{\langle s,s' \rangle} \in \Upsilon^{\eta}},
\ \lambda^{\eta},
\ \{\rho^{\eta}_{s}\}_{s\in S^{\eta}}
\rangle
\]
where
\begin{itemize}

\item $\underline{\eta} = \langle r^{\eta},
\ S^{\eta},
\ \pi^{\eta},
\ \epsilon^{\eta},
\ \Upsilon^{\eta},
\ \{\sigma^{\eta}_{\langle s,s' \rangle}\}_{{\langle s,s' \rangle} \in \Upsilon^{\eta}}
\rangle$ is an $n$-cell shell,
called the {\em base $n$-cell shell} of $\eta$;

\item {\em assignment of cells:}
$\lambda^{\eta} = \coprod_{0\leq i\leq n} \lambda^{\eta}_i$,
where $\lambda^{\eta}_i$ is a function from $S^{\eta}_i$ to $\Sigma_i$
such that for any $s\in S_i$, $\lambda_i(s) \in \Sigma_{i,\epsilon(s)}$

\item {\em identification in boundaries:}
for $s \in S^{\eta}_{i}$,
$\rho^{\eta}_{s}$ is an $i-1$-frame isomorphism
from $\eta|_s$ to $\partial_{i}(\lambda^{\eta}_{i}(s))$

\end{itemize}
which satisfy the following conditions:
\begin{itemize}

\item {\em mutuality:} $\eta|_{r^{\eta}} = 
\langle 
\ S^{\eta|_{r^{\eta}}},
\ \pi^{\eta|_{r^{\eta}}},
\ \epsilon^{\eta|_{r^{\eta}}},
\ \Upsilon^{\eta},
\ \{\sigma^{\eta}_{\langle s,s' \rangle}\}_{{\langle s,s' \rangle} \in \Upsilon^{\eta}},
\ \lambda^{\eta|_{r^{\eta}}},
\ \{\rho^{\eta}_{s}\}_{s\in S^{\eta|_{r^{\eta}}}}
\rangle$
is an $n-1$-frame,
where
  \begin{itemize}

  \item $\langle 
\ S^{\eta|_{r^{\eta}}},
\ \pi^{\eta|_{r^{\eta}}},
\ \epsilon^{\eta|_{r^{\eta}}},
\ \Upsilon^{\eta},
\ \{\sigma^{\eta}_{\langle s,s' \rangle}\}_{{\langle s,s' \rangle} \in \Upsilon^{\eta}} \rangle = \underline{\eta}|_{r^{\eta}}$

  \item $\lambda^{\eta|_{r^{\eta}}} = \lambda^{\eta}|_{S^{\eta|_{r^{\eta}}}}$.

  \end{itemize}


\end{itemize}
\end{Definition}

\begin{Definition}[$n$-frame]
{\em $n$-frame diagram} or {\em $n$-frame} $\zeta$ is 
\[
\langle 
  S^{\zeta},
\ \pi^{\zeta},
\ \epsilon^{\zeta},
\ \Upsilon^{\zeta},
\ \{\sigma^{\zeta}_{\langle s,s' \rangle}\}_{{\langle s,s' \rangle} \in \Upsilon^{\zeta}},
\ \lambda^{\zeta},
\ \{\rho^{\zeta}_{s}\}_{s\in S^{\zeta}}
\rangle
\]
where 
\begin{itemize}

\item $\underline{\zeta} = \langle 
  S^{\zeta},
\ \pi^{\zeta},
\ \epsilon^{\zeta},
\ \Upsilon^{\zeta},
\ \{\sigma^{\zeta}_{\langle s,s' \rangle}\}_{{\langle s,s' \rangle} \in \Upsilon^{\zeta}}
\rangle$ is an $n$-frame shell,
called the {\em base $n$-frame shell} of $\zeta$;

\item {\em assignment of cells:}
$\lambda^{\zeta} = \coprod_{0\leq i\leq n} \lambda^{\zeta}_i$,
where $\lambda^{\zeta}_i$ is a function from $S^{\zeta}_i$ to $\Sigma_i$
such that for any $s\in S_i$, $\lambda_i(s) \in \Sigma_{i,\epsilon(s)}$

\item {\em identification in boundaries:}
for $s \in S^{\zeta}_{i}$,
$\rho^{\zeta}_{s}$ is an $i-1$-frame isomorphism
from $\zeta|_s$ to $\partial_{i}(\lambda^{\zeta}_{i}(s))$

\end{itemize}
which satisfy the following conditions:
\begin{itemize}

\item {\em mutuality:} for every $s\in S^{\zeta}_n$,
$\zeta|^s = \langle
  s,
\ S^{\zeta|^s},
\ \pi^{\zeta|^s},
\ \epsilon^{\zeta|^s},
\ \Upsilon^{\zeta|^s},
\ \{\sigma^{\zeta}_{\langle s,s' \rangle}\}_{{\langle s,s' \rangle} \in \Upsilon^{\zeta|^s}},
\ \lambda^{\zeta|^s},
\ \{\rho^{\zeta}_{t}\}_{t\in S^{\zeta|^s}}
\rangle$ is an $n$-cell diagram, where
  \begin{itemize}

  \item $\langle
  s,
\ S^{\zeta|^s},
\ \pi^{\zeta|^s},
\ \epsilon^{\zeta|^s},
\ \Upsilon^{\zeta|^s},
\ \{\sigma^{\zeta}_{\langle s,s' \rangle}\}_{{\langle s,s' \rangle} \in \Upsilon^{\zeta|^s}} \rangle = \underline{\zeta}|^s$

  \item $\lambda^{\zeta|^s} = \lambda^{\zeta}|_{S^{\zeta|^s}}$

  \end{itemize}


\item {\em compatibility on links:}
      \begin{itemize}
       \item for $\langle s, s' \rangle \in \Upsilon^{\zeta}_{n-1}$,
	     $\lambda(s) = (\lambda(s'))^*$ and
       \item for $\langle s, s' \rangle \in \Upsilon^{\zeta}_{n-1}$
	     and $t \in S^{\zeta|_s}$,
	     $\rho^{\zeta}_{s}(t) = \rho^{\zeta}_{s'}
	     (\sigma_{\langle s,s' \rangle} (t))$
      \end{itemize}
\end{itemize}
\end{Definition}

\begin{Proposition}
For $\langle s, s' \rangle \in \Upsilon^{\zeta}_{n-1}$
and $t \in S^{\zeta}_k$ for some $k\leq n-1$
such that $(\pi^{\zeta})^{n-k-1}(t)=s$,
$(\lambda(t))^* = \lambda(\sigma^{\zeta}_{\langle s, s' \rangle}(t))$.
\end{Proposition}


\begin{Proof}
It is induced from the compatibility on links and
the definition of $(-)^*$ for $n-2$-frames.
\end{Proof}

\begin{Definition}[$\cong_n$, $\Frm{n}$, $(-)^*$]
For two $n$-frames
\begin{align*}
\zeta &=
\langle 
  S^{\zeta},
\ \pi^{\zeta},
\ \epsilon^{\zeta},
\ \Upsilon^{\zeta},
\ \{\sigma^{\zeta}_{\langle s,s' \rangle}\}_{{\langle s,s' \rangle} \in \Upsilon^{\zeta}},
\ \lambda^{\zeta},
\ \{\rho^{\zeta}_{s}\}_{s\in S^{\zeta}}
\rangle\quad\text{and}  \\
\zeta' &=
\langle 
  S^{\zeta'},
\ \pi^{\zeta'},
\ \epsilon^{\zeta'},
\ \Upsilon^{\zeta'},
\ \{\sigma^{\zeta'}_{\langle s,s' \rangle}\}_{{\langle s,s' \rangle} \in \Upsilon^{\zeta'}},
\ \lambda^{\zeta'},
\ \{\rho^{\zeta'}_{s}\}_{s\in S^{\zeta'}}
\rangle\,\text{,}
\end{align*}
an {\em isomorphism} $f$
from $\zeta$ to $\zeta'$ is 
an isomorphism of $n$-frame shells
$f:\underline{\zeta} \rightarrow \underline{\zeta'}$
(with its inverse $f^{-1}$) such that
\begin{itemize}

\item for $s\in S^{\zeta}$,
$\lambda^{\zeta}(s) = \lambda^{\zeta'}(f(s))$
\quad ($\Leftrightarrow$
$\lambda^{\zeta'}(s') = \lambda^{\zeta}(f^{-1}(s'))$),

\item for $s\in S^{\zeta}$ and $t\in S^{\zeta|_s}$,
$\rho^{\zeta}_{s}(t) = \rho^{\zeta'}_{f(s)}(f(t))$
\quad ($\Leftrightarrow$
$\rho^{\zeta'}_{s'}(t') = \rho^{\zeta}_{f^{-1}(s')}(f^{-1}(t'))$).

\end{itemize}
When an isomorphism $f$ from $\zeta$ to $\zeta'$ exists,
we say that $\zeta$ is {\em isomorphic} to $\zeta'$, and write
$f:\zeta \cong_{n} \zeta'$,
$\zeta \cong_{n} \zeta'$, or simply $\zeta \cong \zeta'$.
Obviously $\cong_{n}$ is an equivalence relation.
%
The collection of all $n$-frames is denoted by $\Frm{n}$.
For an $n$-frame $\zeta$, $(\zeta)^*$ is defined as
$\langle 
  S^{\zeta},
\ \pi^{\zeta},
\ \epsilon^{(\zeta)^*},
\ \Upsilon^{\zeta},
\ \{\sigma^{\zeta}_{\langle s,s' \rangle}\}_{{\langle s,s' \rangle} \in \Upsilon^{\zeta}},
\ \lambda^{(\zeta)^*}
\rangle$ where $\epsilon^{(\zeta)^*}(s) = - \epsilon^{\zeta}(s)$ and
$\lambda^{(\zeta)^*}(s) = (\lambda^{\zeta}(s))^*$.
It is easy to check well-definedness,
that is, $(\zeta)^*$ is in fact an $n$-frame, and $((\zeta)^*)^*=\zeta$.
\end{Definition}
An $n$-cell diagram can be seen as an $n$-frame.
Thus we can define an {\em isomorphism} between $n$-cell diagrams similarly.

\begin{Remark}
Indeed, conditions for $\rho$ in the definitions of
cell diagrams and frames {\em ensure} the commutativity of links
and other commutativity of their base shells (it is easy to check this).
Therefore if we use shells only for diagrams,
we need not introduce such commutativity.
A main purpose to do it is
to treat closure operations for shells.
Due to commutativity, a closure becomes unique in a sense.
\end{Remark}

\subsection{$\omega$-hypergraphs}

\begin{Definition}[$\omega$-hypergraph]
An {\em $\omega$-hypergraph} $G = \langle \Sigma, \partial \rangle$
consists of
\begin{itemize}

\item $\Sigma = \coprod_{0\leq i} \Sigma_{i}$, and

\item $\partial = \coprod_{1\leq i} \partial_{i}$.

\end{itemize}

\end{Definition}

\begin{Remark}
Boundaries $\partial_{i}$ depend on frames in the previous step of
the inductive definition.  Therefore as pointed out in \cite{HM},
to formalize the definition of $\omega$-hypergraphs in a logical system,
we need a sort of {\em dependent choice} axiom, $\textit{GDC}_\tau$ in 
\cite{SF} \S 4.4.3 or $\textit{DC}_1$ in \cite{AF} \S 8.2.3.
The strength of this is in between the countable axiom of choice
and the full axiom of choice .
\end{Remark}

\section{Pasting diagrams and their closures}

\begin{Definition}[$n$-pasting shells]
An {\em $n$-pasting shell} consists of the same data
and conditions as an $n$-frame shell, but at the last induction step,
the closedness condition is not required.
That is, $n-1$-nodes which do not appear
in $\Upsilon^{\xi}_{n-1}$ are allowed. 
We call them {\em open} nodes of the pasting shell.
An $n$-pasting shell is {\em positive} or {\em negative}
if for all $s\in S^{\xi}_n$, $\epsilon^{\xi}_n(s)=1$ or $-1$, respectively.
\end{Definition}

\begin{Definition}[$n$-pasting diagrams]
An {\em $n$-pasting diagram} $\zeta$ is defined
in the same way as $n$-frame,
but $\underline{\zeta}$ is an $n$-pasting shell
instead of an $n$-frame shell.
$\cong_n$, $\PD{n}$, $(-)^*$ is also defined similarly.
An $n$-pasting diagram is {\em positive} or {\em negative}
if the base $n$-pasting shell is positive or negative, respectively.
\end{Definition}

\begin{Lemma}\label{lemma:closure}
Consider an $n$-pasting shell
$\xi=
\langle 
  S^{\xi},
\ \pi^{\xi},
\ \epsilon^{\xi},
\ \Upsilon^{\xi},
\ \{\sigma^{\xi}_{\langle s,s' \rangle}\}_{{\langle s,s' \rangle} \in \Upsilon^{\xi}}
\rangle$.
Let a condition $\Psi(y_0, y_1,\ldots, y_m; x_0, x_1,\ldots, x_m)$
($1\leq m$) be abbreviated that 
\begin{itemize}

\item $\langle x_{0},x_{1} \rangle,
\langle x_{2},x_{3} \rangle,\ldots,
\langle x_{m-1},x_{m} \rangle \in \Upsilon^{\xi}_{n-2}$,

\item $\langle y_{0},y_{1} \rangle,
\langle y_{2},y_{3} \rangle,\ldots,
\langle y_{m-1},y_{m} \rangle \in \Upsilon^{\xi}_{n-1}$,

\item $\pi^{\xi}(x_{i})=y_{i}$ and

\item $\sigma^{\xi}_{\langle x_{m-1},x_{m} \rangle} \circ\cdots\circ
\sigma^{\xi}_{\langle x_{2},x_{3} \rangle} \circ
\sigma^{\xi}_{\langle y_{1},y_{2} \rangle} \circ
\sigma^{\xi}_{\langle x_{0},x_{1} \rangle}(x_{0}) = x_{m}$.


\end{itemize}
Note that same nodes may be duplicated in parameters of $\Psi$;
in paticular, $y_0$ may be equal to $y_m$. Then
\begin{enumerate}

\item For every $y_0, y_1,\ldots, y_m$ and $x_0, x_1,\ldots, x_m$ satisfying
$\Psi(y_0, y_1,\ldots, y_m; x_0, x_1,\ldots, x_m)$, we have
$\Psi(y_m, y_{m-1},\ldots, y_0; x_m, x_{m-1},\ldots, x_0)$.

\item For every open node $y_0$ and its child $x_0$ ,
there uniquely exist
$y_0, y_1,\ldots, y_m$ and $x_0, x_1,\ldots, x_m$ satisfying
$\Psi(y_0, y_1,\ldots, y_m; x_0, x_1,\ldots, x_m)$
and that $y_m$ is an open node
($y_1,\ldots, y_{m-1}$ are not open by the second condition of $\Psi$).

\end{enumerate}
 
\end{Lemma}

\begin{Proof}
(1) Trivial from the conditions of frame shells.
(2) Starting from $y_{0}$ and $x_{0}$, we can uniquely fix
a required sequence
$y_{0},x_{0},x_{1},y_{1},y_{2},x_{2},\ldots$ by the following process:
For $x_{2i}$ ($0\leq i$), $x_{2i+1}$ is uniquely determined
by the bijectivity of links; then $\pi^{\xi}(x_{2i+1}) = y_{2i+1}$
and for $y_{2i+1}$, $y_{2i+2}$ is again uniquely
determined by the bijectivity of links;
therefore for $x_{2i+1}$, 
$\sigma_{\langle y_{2i+1}, y_{2i+2} \rangle}(x_{2i+1})=x_{2i+2}$
is also unique.
Next, we show that this process necessarily gets to an open node.
Every link $\langle x, x' \rangle \in \Upsilon^{\xi}_{n-2}$
appears at most once in the process because
for a link to appear twice means the existence of a link
$\langle y, y_0 \rangle \in \Upsilon^{\xi}_{n-1}$ for some $y$,
and this contradicts that $y_0$ is open.
Since $S_{n-2}$ is finite and so is $\Upsilon^{\xi}_{n-2}$,
the process starting from an open node $y_0$
reaches an open node $y_m$ for a finite $m$ and stops there. 
\end{Proof}

\begin{Proposition}[closer and closure of an $n$-pasting shell]
For any $n$-pasting shell
$\xi=
\langle 
  S^{\xi},
\ \pi^{\xi},
\ \epsilon^{\xi},
\ \Upsilon^{\xi},
\ \{\sigma^{\xi}_{\langle s,s' \rangle}\}_{{\langle s,s' \rangle} \in \Upsilon^{\xi}}
\rangle$,
we can construct a {\em closer} of $\xi$, an $n$-cell shell
$\hat{\xi}=
\langle 
r^{\hat{\xi}},
\ S^{\hat{\xi}},
\ \pi^{\hat{\xi}},
\ \epsilon^{\hat{\xi}},
\ \Upsilon^{\hat{\xi}},
\ \{\sigma^{\hat{\xi}}_{\langle s,s' \rangle}\}_{{\langle s,s' \rangle} \in \Upsilon^{\hat{\xi}}}
\rangle$, and
a {\em closure} of $\xi$, an $n$-frame shell
$\overline{\xi}=
\langle 
  S^{\overline{\xi}},
\ \pi^{\overline{\xi}},
\ \epsilon^{\overline{\xi}},
\ \Upsilon^{\overline{\xi}},
\ \{\sigma^{\overline{\xi}}_{\langle s,s' \rangle}\}_{{\langle s,s' \rangle} \in \Upsilon^{\overline{\xi}}}
\rangle$
uniquely up to isomorphisms and polarity as follows:
Let $\{s_0,s_1,\ldots,s_k\}$ be the set of open nodes of $\xi$.
For each $s_l$, 
we prepare an $n-1$-cell shell
$\tau_{l}=
\langle 
t_l,
\ S^{\tau_{l}},
\ \pi^{\tau_{l}},
\ \epsilon^{\tau_{l}},
\ \Upsilon^{\tau_{l}},
\ \{\sigma^{\tau_{l}}_{\langle t,t' \rangle}\}_{{\langle t,t' \rangle} \in \Upsilon^{\tau_{l}}}
\rangle$,
isomorphic to $(\xi|^{s_l})^*$ via an isomorphsim
$f_l: \tau_{l}\rightarrow (\xi|^{s_l})^*$.
Then the components of the closer $\hat{\xi}$ are:
\begin{itemize}

\item $S^{\hat{\xi}}_n
	= \{r^{\hat{\xi}}\}$ where $\{r^{\hat{\xi}}\}$ is a singleton set, and
      $S^{\hat{\xi}}_i
	= \coprod_{0\leq l\leq k} S^{\tau_{l}}_i$
	for $0\leq i\leq n-1$,

\item $\pi^{\hat{\xi}}_{n-1}
	= \pi_{r^{\hat{\xi}}}$
        where $\pi_{r^{\hat{\xi}}}(t_l) = r^{\hat{\xi}}$ for $0\leq l\leq k$, and
      $\pi^{\hat{\xi}}_i
	= \coprod_{0\leq l\leq k} \pi^{\tau_{l}}_i$
	for $0\leq i\leq n-2$,

\item $\epsilon^{\hat{\xi}}
	= (\coprod_{0\leq l\leq k} \epsilon^{\tau_{l}})
	\amalg \epsilon_{r^{\hat{\xi}}}$ where
      $\epsilon_{r^{\hat{\xi}}}(r^{\hat{\xi}})=-1$
      (the {\em negative} closer) or $1$ (the {\em positive} closer),

\item $\Upsilon^{\hat{\xi}} =
(\coprod_{0\leq l\leq k} \Upsilon^{\tau_{l}}) \amalg
\{\langle f_l^{-1}(x),f_{l'}^{-1}(x') \rangle\,|\,\Phi(s_l, s_{l'}; x,x')\}$,
where $\Phi(s_l, s_{l'}; x,x')$ is abbreviated that
there exist
$y_{0}=s_l, y_{1}, y_{2},\ldots, y_{m}=s_{l'}$ and
$x_{0}=x, x_{1}, x_{2},\ldots, x_{m}=x'$
satisfying $\Psi(y_0, y_1,\ldots, y_m; x_0, x_1,\ldots, x_m)$,

\item $\{\sigma^{\hat{\xi}}_{\langle s,s' \rangle}\}_{{\langle s,s' \rangle} \in \Upsilon^{\hat{\xi}}}$ is defined as
  \begin{itemize}

  \item $\sigma^{\hat{\xi}}_{\langle t,t' \rangle} = 
\sigma^{\tau_l}_{\langle t,t' \rangle}$ for $\langle t,t' \rangle
\in \Upsilon^{\tau_l}$,

  \item $\sigma^{\hat{\xi}}_{\langle f_l^{-1}(x),f_{l'}^{-1}(x')\rangle}=
f_{l'}^{-1} \circ
\sigma^{\xi}_{\langle x_{m},x \rangle} \circ\cdots\circ
\sigma^{\xi}_{\langle x_{2},x_{3} \rangle} \circ
\sigma^{\xi}_{\langle y_{1},y_{2} \rangle} \circ
\sigma^{\xi}_{\langle x_{0},x_{1} \rangle} \circ
f_l$.

  \end{itemize}

\end{itemize}
and the components of the closure $\overline{\xi}$ are:
\begin{itemize}

\item $S^{\overline{\xi}} = S^{\xi} \amalg S^{\hat{\xi}}$,

\item $\pi^{\overline{\xi}} = \pi^{\xi} \amalg \pi^{\hat{\xi}}$,

\item $\epsilon^{\overline{\xi}} = \epsilon^{\xi} \amalg \epsilon^{\hat{\xi}}$,

\item $\Upsilon^{\overline{\xi}} = \Upsilon^{\xi} \amalg
\{\langle s_l,t_l \rangle, \langle t_l,s_l \rangle\,|\,0\leq l\leq k\} \amalg
\Upsilon^{\hat{\xi}}$,

\item $\{\sigma^{\overline{\xi}}_{\langle s,s' \rangle}\}_{{\langle s,s' \rangle} \in \Upsilon^{\overline{\xi}}}$ is defined as
  \begin{itemize}

  \item $\sigma^{\overline{\xi}}_{\langle s,s' \rangle} = 
\sigma^{\xi}_{\langle s,s' \rangle}$ for $\langle s,s' \rangle
\in \Upsilon^{\xi}$,

  \item $\sigma^{\overline{\xi}}_{\langle t_l,s_l \rangle}=f_l$ and
$\sigma^{\overline{\xi}}_{\langle s_l,t_l \rangle}=f_l^{-1}$,

  \item $\sigma^{\overline{\xi}}_{\langle t,t' \rangle} = 
\sigma^{\hat{\xi}}_{\langle t,t' \rangle}$
for $\langle t,t' \rangle \in \Upsilon^{\hat{\xi}}$.

  \end{itemize}
\end{itemize}

\end{Proposition}


\begin{Proof}
First, we will check conditions for 
$\langle f_l^{-1}(x),f_{l'}^{-1}(x') \rangle$ and 
$\sigma^{\overline{\xi}}_{\langle f_l^{-1}(x),f_{l'}^{-1}(x')\rangle}$.
Other parts are rather easy:
\begin{itemize}

\item The mutuality condition is obvious.
The bijectivity condition is shown by the Lemma \ref{lemma:closure} (2) and
the involution condition by the Lemma \ref{lemma:closure} (1).

\item The conjugation condition is derived from the following results:
$\epsilon^{\xi}(x_{2i})\epsilon^{\xi}(x_{2i+1})=-1$
by the conjugation of $\xi$,
$\epsilon^{\xi}(x_{2i+1})\epsilon^{\xi}(x_{2i+2})=-1$
by the definition of $\sigma^{\xi}_{\langle y_{2i+1},y_{2i+2} \rangle}$,
$\epsilon^{\tau_l}(f_l^{-1}(x))\epsilon^{\xi}(x)=-1$ and
$\epsilon^{\tau_{l'}}(f_{l'}^{-1}(x'))\epsilon^{\xi}(x')=-1$.

\item The correspondence of links condition
for $\langle f_l^{-1}(x),f_{l'}^{-1}(x') \rangle$
is shown by chaining the correspondence of links in $\xi$ and 
commutativity of $f$, $f'$.

\item The closedness condition is straightforward from the Lemma
      \ref{lemma:closure} (2).

\end{itemize}


Next, we will check the commutativity of links condition for the closure.
The following three cases are possible:
\begin{enumerate}

\item All links are in $\Upsilon^{\xi}$;\label{comm:case1}

\item All links are in $\Upsilon^{\hat{\xi}}$;\label{comm:case2}

\item Links in
$\{\langle s_l,t_l \rangle, \langle t_l,s_l \rangle\,|\,0\leq l\leq k\}$
occur.\label{comm:case3}

\end{enumerate}
 \begin{figure}
  \begin{center}
  \includegraphics{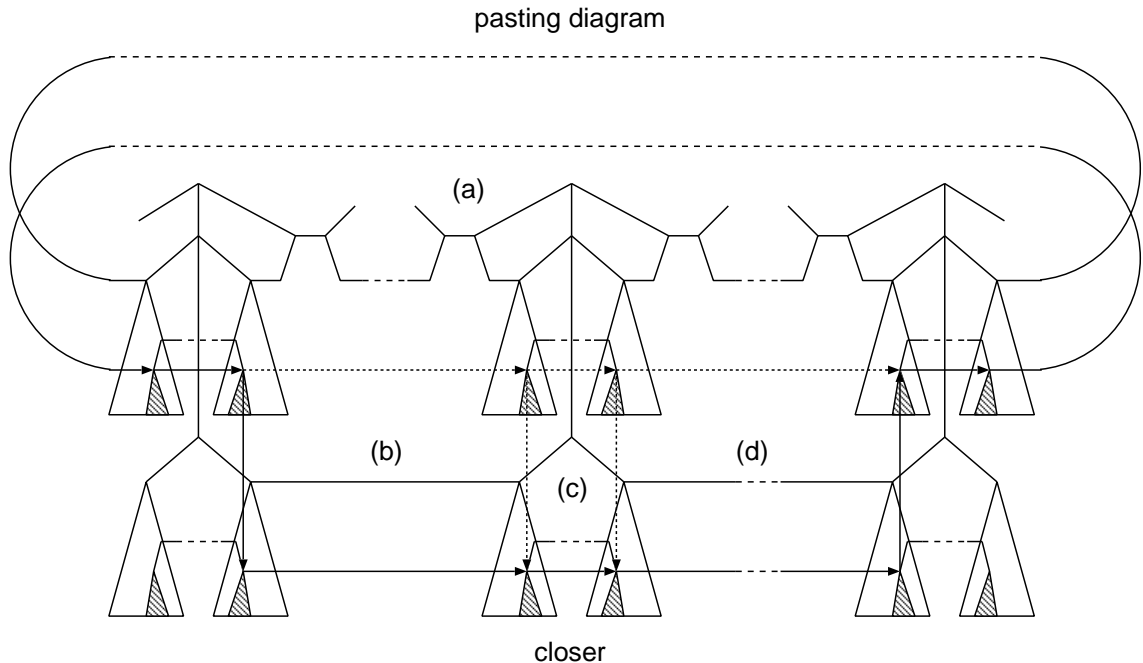}
  \end{center}
  \caption{A chain of shell isomorphisms in a closure\label{comm:fig1}}
 \end{figure}
Commutativity for the case \ref{comm:case1} is trivial
from the definition of pasting diagrams.
For the cases \ref{comm:case2} and \ref{comm:case3},
we show commutativity of a path in Figure \ref{comm:fig1}.
The oval path (a) is commutative from the definition of pasting diagrams;
the square (b) and (c) is from the construction of closers;
and the square (d) is by pasting an alternation of
the type (b) and (c) squares.
Thus the outer path of arrows in this case is commutative by pasting (a)--(d).
In the general case \ref{comm:case3}, 
the side trip (b)--(d) might occur several times.
The largest roundabouts are paths running only through the closer.
This implies commutativity for the case \ref{comm:case2}, that is,
the commutativity of links condition for the closer.
\end{Proof}

\section{Directed $\omega$-hypergraphs}

\subsection{directed $\omega$-hypergraphs}

\begin{figure}
 \begin{center}
  \includegraphics{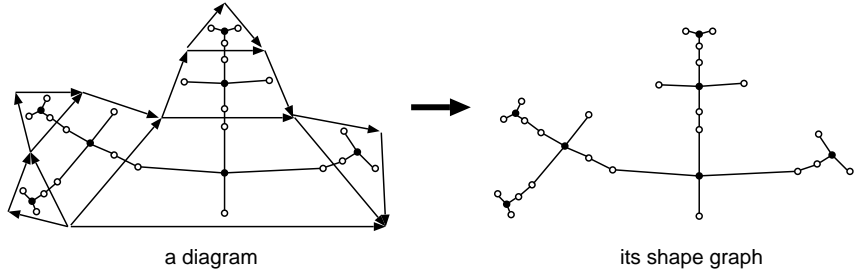}
 \end{center}
 \caption{Extraction of shape graphs\label{shape:fig1}}
\end{figure}
\begin{Definition}[shape graph]
The (undirected) {\em shape graph} $\langle N_0^{\xi}, N_1^{\xi}
E_0^{\xi}, E_1^{\xi} \rangle$
of an $n$-frame $\xi$ is defined as follows:
\begin{itemize}

\item the body node set $N_0^{\xi}$ is $S^{\xi}_n$;

\item the foot node set $N_1^{\xi}$ is $S^{\xi}_{n-1}$;

\item the leg edge set $E_0^{\xi}$ is
$\{\{s,t\}\,|\, s\in S^{\xi}_{n}, t\in S^{\xi}_{n-1}, \pi^{\xi}(t)=s \}$;

\item the link edge set $E_1^{\xi}$ is
$\{\{t,t'\}\,|\, \langle t,t'\rangle \in \Upsilon_{n-1}^{\xi}\}$
(note that $\pi^{\xi}(t) \neq \pi^{\xi}(t')$
from the definition of $\Upsilon^{\xi}_{n-1}$).

\end{itemize}
If $s\in N_0^{\xi}$, $t\in N_1^{\xi}$ and $\{s,t\}\in E_0^{\xi}$
then we call $t$ a {\em foot} of $s$ and $\{s,t\}$ a {\em leg} of $s$.
The {\em shape graph} $\langle N_0^{\zeta}, N_1^{\zeta}
E_0^{\zeta}, E_1^{\zeta} \rangle$ of an $n$-frame $\zeta$ is the shape
graph of $\underline{\zeta}$.

The {\em shape graph} of $n$-cell shells or $n$-cell diagrams is defined 
as a special case of $n$-frame shells or $n$-frames.
\end{Definition}

\begin{Definition}[$n$-directed $n$-frame]
An $n$-frame $\zeta$ is {\em $n$-directed}
if it satisfies the following conditions:
\begin{itemize}

\item {\em headedness:}
for exactly one $s\in S^{\zeta}_n$
and every other $s'\in S^{\zeta}_n$, either
$\lambda^{\zeta}_n(s)$ is positive and 
$\lambda^{\zeta}_n(s')$ is negative or
$\lambda^{\zeta}_n(s)$ is negative and 
$\lambda^{\zeta}_n(s')$ is positive,
where $s$ is called the {\em positive} or {\em negative head}
of $\zeta$, respectively;

\item {\em connectedness:}
its shape graph is connected;

\item {\em acyclicity:}
the graph obtained from its shape graph
by getting rid of a body node corresponding to the head, its legs and feet, 
and link edges connected to them, is acyclic
(indeed, this graph is a tree).

\end{itemize}
An $n$-frame with the positive head is 
said to be {\em positively $n$-directed},
and that with the negative head be {\em negatively $n$-directed}.
The same $n$-frame can be both positively and negatively $n$-directed.
An $n+1$-cell whose boundary is
such an $n$-frame is called a {\em simple} $n+1$-cell.
\end{Definition}

\begin{figure}
 \begin{center}
  \includegraphics{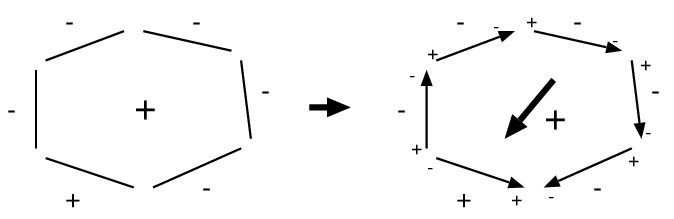}
 \end{center}
 \caption{Direction via polarity \label{directed:fig1}}
\end{figure}

\begin{Definition}[directed $n$- and $\omega$-hypergraph]
A {\em directed $\omega$-hypergraph} is an $\omega$-hypergraph
which satisfies the following condition:
\begin{itemize}

\item {\em directedness:} for each $i\geq 1$,
the boundary of any positive $i$-cell is 
a positively $i-1$-directed $i-1$-frame,
and that of any negative $i$-cell is
a negatively $i-1$-directed $i-1$-frame.

\end{itemize}
For each $n\geq 0$, a {\em directed $n$-hypergraph} is also
defined as an $n$-hypergraph satifying the same condition.
\end{Definition}

\subsection{directed shells}

In the category theory besed on $n$- or $\omega$-hypergraphs,
{\em directed $n$-cell shells} and {\em directed $n$-frame shells}
play the role of shape diagrams in the usual theory.
They are defined by adding some conditions to the induction step of
the definitions of $n$-cell shells and $n$-frame shells

\begin{Definition}[directed $n$-cell shell]
An additional condition is as follows:
\begin{itemize}

\item {\em directedness:} If $\epsilon(r)=1$, then $\theta|_r$
is a positively directed $n-1$-frame shell and
if $\epsilon(r)=-1$, then it is a negatively directed one.

\end{itemize}
\end{Definition}

\begin{Definition}[directed $n$-frame shell]
Additional conditions are as follows:
\begin{itemize}

\item {\em headedness:}
for exactly one $s\in S^{\xi}_n$
and every other $s'\in S^{\xi}_n$, either
$\epsilon^{\xi}_n(s)=1$ and 
$\epsilon^{\xi}_n(s')=-1$ or
$\epsilon^{\xi}_n(s)=-1$ and 
$\epsilon^{\xi}_n(s')=1$,
where $s$ is called the {\em positive} or {\em negative head}
of $\xi$, respectively;

\item {\em connectedness:} its shape graph is connected;

\item {\em acyclicity:}
the graph obtained from its shape graph
by getting rid of a body node corresponding to the head, its legs and feet, 
and link edges connected to them, is acyclic
(indeed, this graph is a tree).

\end{itemize}
An $n$-frame shell with the positive head is 
said to be {\em positively directed},
and that with the negative head be {\em negatively directed}.
The same $n$-frame shell can be both positively and negatively directed.
If for an $n+1$-cell shell $\theta$ with root $r$,
$\theta|_r$ is such an $n$-frame shell,
then it is called a {\em simple} $n+1$-cell shell.
\end{Definition}

\begin{Proposition}
An $n$- or $\omega$-hypergraph $\langle \Sigma, \partial \rangle$
is a directed $n$- or $\omega$-hypergraph iff
for each positive $i$-cell $c$, $\underline{\partial_i(c)}$ is
a positivery directed $i-1$-frame shell
and for each negative $i$-cell, it is a negatively directed one.
\end{Proposition}

\begin{Proof} By induction on dimensions.
\end{Proof}

\begin{Definition}[directed $n$-frame and directed $n$-cell diagram]
An $n$-frame $\zeta$ is 
a {\em positively} or {\em negatively directed} $n$-frame
if $\underline{\zeta}$ is
a positively or negatively directed $n$-frame shell, respectively.
Also a {\em positively} or {\em negatively directed}
$n$-cell diagram is defined in the same way.
\end{Definition}

\begin{Corollary}
In any directed $n$- or $\omega$-hypergraph,
an $n$-frame is a positively or negatively directed $n$-frame iff
it is an positively or negatively $n$-directed $n$-frame, respectively.
\end{Corollary}

\begin{Definition}[$\DFrm{k}$]
For a directed $n$- or $\omega$-hypergraph,
the category (groupoid) whose objects are all directed $k$-frames and
whose arrows are all isomorphisms is denoted by $\DFrm{k}$.
The collection of all directed $k$-frames is also denoted by $\DFrm{k}$.
\end{Definition}

In the rest of this paper, we will mainly use the usual diagramatic
notations for shells and diagrams (Figure \ref{directed:fig1}).

\subsection{examples}

\begin{Example}[hypergraph in rewriting]
An (directed) {\em hypergraph} used in hypergraph rewritinng \cite{DP}
is a directed $1$-hypergraph.
\end{Example}

\begin{Example}[$\omega$-multigraph]
An {\em $\omega$-multigraph} is a directed $\omega$-hypergraph
$\langle \Sigma, \partial \rangle$ such that
$\Sigma_{0,1}$ is a singlton set and that
any $c\in\Sigma_{1,1}$ is a simple $1$-cell.
\end{Example}

\begin{Example}[doublegraph]
Doublegraphs are underlying graph-like structures for double categories.
They are obtained by splitting $1$-cells into
vertical cells and horizontal cells:
\begin{gather*}
\Sigma_{0} = \Sigma_{0,1} \amalg \Sigma_{0,-1} \\
\Sigma_{1} = \Sigma_{1,1} \amalg \Sigma_{1,-1} \\
\Sigma_{1,1} = \Sigma_{1,1}^{v} \amalg \Sigma_{1,1}^{h} \\
\Sigma_{1,-1} = \Sigma_{1,-1}^{v} \amalg \Sigma_{1,-1}^{h}
\end{gather*}
and if $c\in \Sigma_{1,k}$, then $c^\ast \in \Sigma_{1,-k}$, etc.
$2$-cells are as follows:
\[
\includegraphics{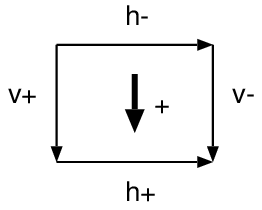}
\]
This appoach can be easily extended to multiple categories.
\end{Example}

\begin{Example}[\textbf{fc}-multigraph]
\textbf{fc}-multigraphs are underlying graph-like structures
for \textbf{fc}-multicategories introduced by T. Leinster \cite{TL}.
Similar notions also appear in \cite{CH}.
They are a mixture of $2$-multigraphs and double graphs.
0-cells and 1-cells are the same as double graphs.
$2$-cells are as follows
\[
\includegraphics{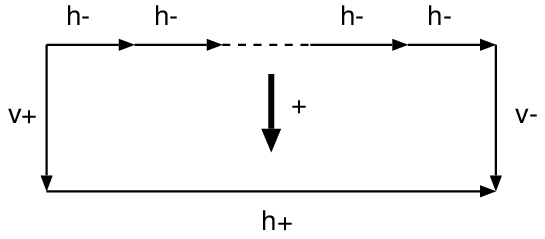}
\]
\end{Example}

\subsection{directed pasting shells and diagrams}

\begin{Definition}[boundary graphs]
The {\em shape graph} of an $n$-pasting shell $\xi$ is defined 
in the same way of $n$-frame shells.
The {\em boundary graph} of an $n$-pasting shell $\xi$ is defined as
the shape graph of $(\hat{\xi}|_{r^{\hat{\xi}}})^\ast$.
\end{Definition}
Note that the {\em boundary graph} of an $n$-cell shell $\theta$ 
as a special case of $n$-pasting shells, matches with
the shape graph of $\theta|_{r^{\theta}}$.

\begin{Definition}[directed $n$-pasting shells]
A {\em directed $n$-pasting shell} $\xi$ is 
an $n$-pasting shell consisting of directed $n$-cell shells
satisfying the following conditions:
\begin{itemize}


\item {\em homegeneity:} it is positive or negative as an $n$-pasting shell.

\item {\em connectedness:} its shape graph is connected;


\item {\em acyclicity:} its shape graph is acyclic
(indeed, this graph is a tree).

\end{itemize}
\end{Definition}

\begin{Definition}[directed $n$-pasting diagrams]
A {\em directed $n$-pasting diagram} is an $n$-pasting diagram
whose base $n$-pasting shell is a directed $n$-pasting shell.
A {\em positive} or {\em negative} directed $n$-pasting diagram
is trivially defined, respectively.
\end{Definition}

\begin{Proposition}
For any directed $n$-frame $\zeta$ and its head node $h$,
we can uniquely split it into an directed $n$-cell diagram $\cod(\zeta)$,
called the {\em codomain} of $\zeta$,
and an $n$-pasting diagram $\dom(\zeta)$,
called the {\em domain} of $\zeta$,
where
\begin{itemize}

\item $\cod(\zeta)=\zeta|^{h}$

\item $\dom(\zeta)$ is obtained by
deleting, from the data of $\zeta$, $\zeta|^h$ and 
indice $\langle s,s' \rangle$ of which
$s$ or $s'$ is in $\Upsilon^{\zeta|^{h}}$.

\end{itemize}
Of course, if $\zeta$ is positively directed,
then $\cod(\zeta)$ is positive and $\dom(\zeta)$ is negative,
and if negatively directed, then
$\cod(\zeta)$ is negative and $\dom(\zeta)$ is positive.
\end{Proposition}

\begin{figure}
 \begin{center}
  \includegraphics{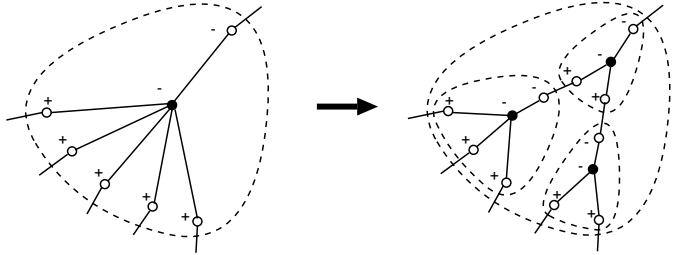}
 \end{center}
 \caption{replacement of the boundary graph\label{replace:fig1}}
\end{figure}
\begin{Proposition}[closure of a directed $n$-pasting shell]
For any negative (resp. positive) directed $n$-pasting shell $\xi$,
(1) its positive (resp. negative) closer $\hat{\xi}$
is a directed $n$-cell shell and
(2) its positive (resp. negative) closure $\overline{\xi}$
is a directed $n$-frame shell.
\end{Proposition}

\begin{Proof}
(1)
We have to check that  $\hat{\xi}|_{r^{\hat{\xi}}}$ is 
a positively directed $n-1$-frame shell.
(i) {\em Headedness:} 
From homogeneity, connectedness and acyclicity of $\xi$,
the set $O_\xi$ of open nodes of $\xi$ contains just one positive node.
Therefore $\hat{\xi}|_{r^{\hat{\xi}}}$ contains just one negative node.
(ii) {\em Connectedness and acyclicity:} 
By induction on the construction of directed $n$-pasting shells.
(a) The boundary graph of a directed $n$-cell shell 
satisfies the connectedness and acyclicity conditions.
(b) If you make a directed $n$-pasting shell $\xi'$ by linking
a directed $n$-pasting shell $\xi$
and a directed $n$-cell shell $\theta$ at one open $n-1$-node
(with satisfying the homogeneity, connectedness and acyclicity for $\xi'$),
then replacement of the boundary graph occurs (Figure \ref{replace:fig1}) and
the resulting boundary graph of $\xi'$ also satisfies
the connectedness and acyclicity.
The acyclicity is obtained by reduction to absurdity.
Suppose the existance of cycles and consider the graph
obtained by deleting link edges of cycles from
the boundary graph, and its polarity of foot nodes, it contradicts
the directedness of the boundary $n-1$-frame shell $\theta|_{r^{\theta}}$
of each $n$-cell shell $\theta$ in the directed $n$-pasting shell.
(c) Any directed $n$-pasting shell can be constructed by finitely
iterating this process, and then the boundary graph of
the closer satisfies the connectedness and acyclicity.
Thus $\hat{\xi}|_{r^{\hat{\xi}}}$ is a positively directed $n-1$-frame shell.
Since $\epsilon^{\hat{\xi}}(r^{\hat{\xi}})=1$, 
$\hat{\xi}$ is a positively directed $n$-cell shell.
The negative case is in parallel.

(2) From the homogeneity, connectedness and acyclicity of $\xi$,
the headedness, connectedness and acyclicity of $\overline{\xi}$ is obvious.
\end{Proof}

\section{Weak $\omega$-categories}

\subsection{$\omega$-identity, $\omega$-invertibility and $\omega$-universality}

We will {\em coinductively} define three notions:
{\em $\omega$-identity}, {\em $\omega$-invertibility} and
{\em $\omega$-universality}.
All $n$-dimensional notions depend on $n+1$- or $n+2$-dimensional ones.
The reader unfamiliar with coinductive definitions may think of
only the case in which coinduction steps terminate. 

One source of our idea is
Michael Makkai's work on anabicategories \cite{MM1}.
At a glance, as Makkai pointed out in \cite{MM2}, 
{\em saturated} anabicategories could be regarded as
2-dimensional weak cateogires of Baez-Dolan.
But we don't think of them to be equivalent notions for some reasons:

1) In anabicategories,
two composite arrows for the same composable sequence of arrows
are equivalent.
While in Baez-Dolan's
there are two opposite universal 2-cells between two composite arrows
by virtue of balancedness, no explicit relation between them appears.
In fact, we can prove that they are equivalences in a sense,
because Baez and Dolan only think of finite dimensional cases.
But we cannot prove it in that way for infinite dimensional cases.
Therefore we need to characterize those opposites as a sort of equivalences.

2) Different from anabicategories,
composites of empty sequence are introduced in Baez-Dolan's
and expected to play the role of identities.
But as well as the above,
we cannot prove the property of identity in infinite dimensional cases.
Thus we also have to define identities explicitly.

3) In (not necessarily saturated) anabicategories,
an object isomorphic to a composite might not be a composite.
It suggests that if we introduce a sort of equivalences,
it is natural to treat equivalences and composition separatedly
and add a saturatedness condition.



\begin{Definition}[$\omega$-identical cells]
An $n$-cell $c$ is {\em $\omega$-identical} if
it is simple and satisfies $\dom(c)\cong\cod(c)$
and the following conditions:
\[
\begin{array}{lcll}
\forall f & & \exists \alpha & \\
{\diagram
& \drtwocell<0>^{f}{\omit} & \\
\urtwocell<0>^{c}{\omit}
\rrtwocell<0>_{f}{\omit} & & \\
\enddiagram} &

\quad\Rightarrow\quad &

{\diagram
& \drtwocell<0>^{f}{\omit} & \\
\urtwocell<0>^{c}{\omit}
\rrtwocell<0>_{f}{<-2.5>\alpha} & & \\
\enddiagram} &
\alpha\text{: $\omega$-universal}
\end{array}
\]
and
\[
\begin{array}{lcll}
\forall g & & \exists \beta & \\
{\diagram
& \drtwocell<0>^{c}{\omit} & \\
\urtwocell<0>^{g}{\omit}
\rrtwocell<0>_{g}{\omit} & & \\
\enddiagram} &

\quad\Rightarrow\quad &

{\diagram
& \drtwocell<0>^{c}{\omit} & \\
\urtwocell<0>^{g}{\omit}
\rrtwocell<0>_{g}{<-2.5>\beta} & & \\
\enddiagram} &
\beta\text{: $\omega$-universal.} \\
\end{array}
\]
\end{Definition}

\begin{Definition}[$\omega$-invertibility, $\omega$-equivalence, $\simeq$]
A pair of $n$-cells $f$ and $g$ is an {\em $\omega$-invertible pair}
if both $f$ and $g$ are simple and satisfy
$\dom(f)\cong\cod(g)$ and $\dom(g)\cong\cod(f)$
and the following condtions:
\[
\begin{array}{ll}
\exists \alpha, i  & \\
{\diagram
& \drtwocell<0>^{f}{\omit} & \\
\urtwocell<0>^{g}{\omit}
\rrtwocell<0>_{i}{<-2.5>\alpha} & & \\
\enddiagram} &
i\text{: $\omega$-identical,} \\
\end{array}
\]
and
\[
\begin{array}{ll}
\exists \beta, j  & \\
{\diagram
& \drtwocell<0>^{g}{\omit} & \\
\urtwocell<0>^{f}{\omit}
\rrtwocell<0>_{j}{<-2.5>\beta} & & \\
\enddiagram} &
j\text{: $\omega$-identical.}
\end{array}
\]
$f$ and $g$ are called {\em $\omega$-invertible}.
We say that 
two $n-1$ cells $\lambda(r^{\dom(f)})$ and $\lambda(r^{\cod(f)})$
are {\em $\omega$-equivalent} and write it as
$\lambda(r^{\dom(f)}) \simeq \lambda(r^{\cod(f)})$.

\end{Definition}

\begin{Definition}[$\omega$-universal cells]
An $n$-cell $u$ is {\em $\omega$-universal} if for any $n$-cell $f$ of
\[
{\diagram
\ddtwocell<0>_{f}{\omit}
\drto^{u} &\\
& \\
& \\
\enddiagram}
\]
there exist an $n$-cell $g$ and an $\omega$-universal $n+1$-cell $\alpha$ to be
\[
{\diagram
\ddtwocell<0>_{f}{<-2>\alpha}
\drto^{u} & \\
& \dlto^g \\
& \\
\enddiagram}
\]
and for two such pairs, $g$ and $\alpha$, $h$ and $\beta$
\[
{\diagram
\ddtwocell<0>_{f}{<-2>\alpha}
\drto^{u} & \\
& \dlto^g \\
& \\
\enddiagram}
\quad
{\diagram
\ddtwocell<0>_{f}{<-2>\beta}
\drto^{u} & \\
& \dlto^h \\
& \\
\enddiagram}
\]
there exist an $\omega$-invertible pair of $n+1$-cells,
$\gamma$ and $\delta$, and two $\omega$-universal $n+2$-cells,
$\Phi$ and $\Psi$, such that
\[
{\diagram
& \\
& \dltwocell^{g}_{h}{\gamma} \\
& \\
\enddiagram}
{\diagram
& \\
& \dltwocell^{h}_{g}{\delta} \\
& \\
\enddiagram} \quad
{\diagram
\ddtwocell<0>_{f}{<-2>\beta}
\drto^{u} & \\
& \dltwocell<0>_{h}{\omit}
\dluppertwocell^{g}{\gamma} \\
& \\
\enddiagram}
{\xy(0,0)
\ar @3 (4,0)^{\Phi}
\endxy}
{\diagram
\ddtwocell<0>_{f}{<-2>\alpha}
\drto^{u} & \\
& \dluppertwocell^{g}{\omit} \\
& \\
\enddiagram} \quad
{\diagram
\ddtwocell<0>_{f}{<-2>\alpha}
\drto^{u} & \\
& \dltwocell<0>_{g}{\omit}
\dluppertwocell^{h}{\delta} \\
& \\
\enddiagram}
{\xy(0,0)
\ar @3 (4,0)^{\Psi}
\endxy}
{\diagram
\ddtwocell<0>_{f}{<-2>\beta}
\drto^{u} & \\
& \dluppertwocell^{h}{\omit} \\
& \\
\enddiagram}
\]
\end{Definition}

\subsection{weak $\omega$-categories}

\begin{Definition}[weak $\omega$-categories]
A directed $\omega$-hypergraph is a {\em weak $\omega$-category}
if it satisfies the following conditions:
\begin{itemize}

\item {\em existence of closers and occupants:}
For any $n$-pasting diagram $P$,
\[
\includegraphics{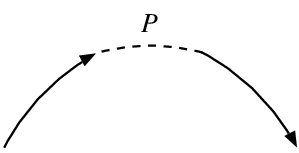} 
\]
there exist
an $\omega$-universal $n+1$-cell $\alpha$ and an $n+1$-cell $h$
such that $\dom(\alpha)\cong_{n}P$ and $\lambda(r^{\cod(\alpha)})=h$:
\[
\includegraphics{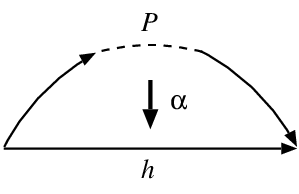} 
\]
Moreover, if $P$ is an empty pasting diagram,
then for each $n-1$-cell diagram $x$,
there exist such $\alpha$ and $h$ satisfying
$\dom(h)\cong_{n-1}\cod(h)\cong_{n-1} x$.
\[
\includegraphics{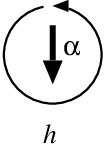} 
\]
We call $\alpha$ an {\em occupant} for $P$ and $h$ a {\em closer} of $P$.

\item {\em weak uniqueness of closers and occupants:}
For two such pairs as above, $\alpha$ and $h$, $\beta$ and $k$
\[
\includegraphics{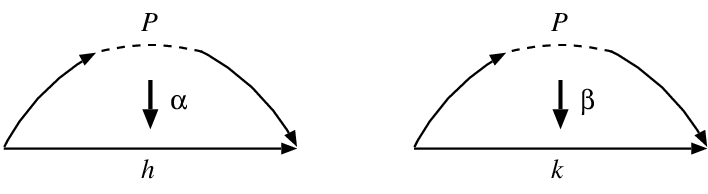} 
\]
there exist an $\omega$-invertible pair  $\gamma$ and $\delta$
and two $\omega$-universal cells $\Phi$ and $\Psi$ such that
\[
\includegraphics{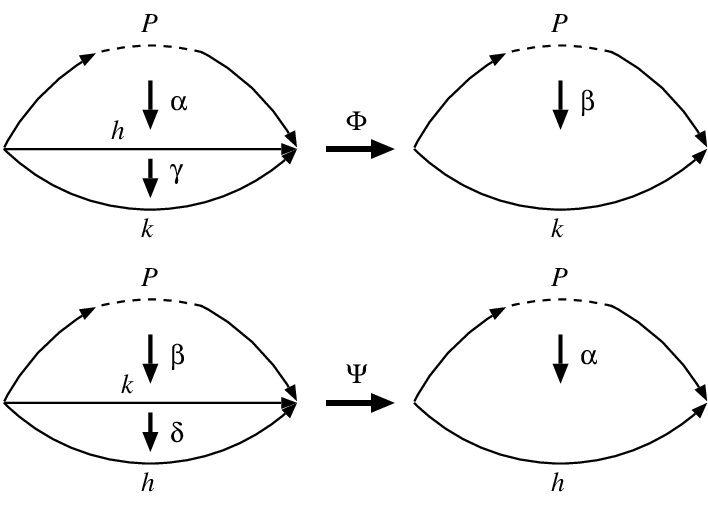} 
\]

\item {\em saturation of closers and occupants:}
For an $\omega$-universal cell $\alpha$ and
an $\omega$-invertible pair $\gamma$ and $\delta$
\[
\includegraphics{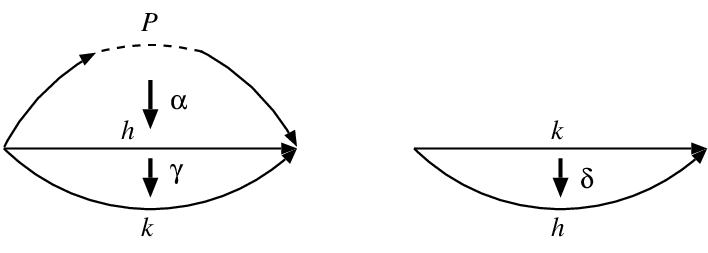} 
\]
there exists $\omega$-universal cells $\beta$, $\Phi$, $\Psi$ such that
\[
\includegraphics{wc5} 
\]

\item {\em $\omega$-identical closers (1):}
For any $\omega$-universal $n+1$-cell $\alpha$ as follows,
an $n$-cell $i$ is $\omega$-identitical:
\[
\includegraphics{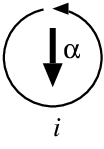}
\]

\item {\em $\omega$-identical closers (2):}
For any $\omega$-identitical $n$-cell $i$,
there is an $\omega$-universal $n+1$-cell $\alpha$ as above.

\item {\em $\omega$-universal closers:}
Any closer for a $n$-pasting diagram made of
$\omega$-universal $n$-cells is $\omega$-universal.

\end{itemize}
A weak $\omega$-categories is defined when
this coinductive definition makes sense.
\end{Definition}

\begin{Proposition}
An $\omega$-identical cell is $\omega$-invertible.
\end{Proposition}

\begin{Proposition}
Every two $\omega$-identical cells are $\omega$-equivalent.
\end{Proposition}

Again, the reader unfamiliar with coinduction may think
of weak $n$-categories.

\begin{Definition}[weak $n$-categories]
A weak $\omega$-category is a {\em weak $n$-category}
if for each $k$ higher than $n$,
all simple $k$-cells are $\omega$-invertible.
\end{Definition}

From the axioms above, we can recognize that
identity is independent from the definition of composition.
In fact, to define our weak $\omega$-categories,
we can exclude empty pasting diagram and related axioms
and add an axiom for existance of $\omega$-identical $n+1$-cells
for each $n$-cells. This is slightly simpler than those defined above.
And we conjecture that our definition would be
equivalent to that of J. Penon \cite{JP},
and furthermore that if we abandon saturatedness and
for each $n$-cell we choose just one $\omega$-identical $n+1$-cell
whose domain and codomain are that $n$-cell,
then the category of our small weak $\omega$-categories 
(with suitable functors) would be isomorphic to
the category of Penon's.

\bigskip
\noindent
\textbf{Acknowledgements:} The first author thanks Prof.\ Sounders
Mac~Lane for his comments and powerful encouragement for his talks at
CT97 and CT99, and Prof.\ Michael Makkai for thoughtful explanation of
his ideas.  We also thank A. Higuchi for helpful discussion.
The first author's work is partially supported by Kyoto Sangyo University.


\begin{thebibliography}{9}

\bibitem{AF} Jeremy Avigad and Solomon Feferman, G\"{o}del's functional
	(``{Dialetica}'') interpretation. In: Samuel R. Buss ed., {\em
	Handbook of Proof Theory} (Elsevier Science B. V., 1998) 337--405.

\bibitem{BD} John C. Baez and James Dolan, Higher-dimensional algebra
	{III}: $n$-categories and the algebra of opetopes, {\em Advances
	in Mathemtaics} 135 (1998) 145--206.

\bibitem{CB} C. Berge, {\em Graphs and Hypergrahs} (North-Holland,
	Amsterdam, 1973).

\bibitem{RF} Ronald Fagin, Degrees of acyclicity for hypergraphs and
	relational database schemes, {\em J. Assoc. Comp. Mach.},
	30(3) (July 1983) 514--550.

\bibitem{SF} Solomon Feferman, Theories of finite type related to
	mathematical practice. In: Jon Barwise ed., {\em Handbook of
	Mathematical Logic} (Elsevier Science B.V., 1977) 913--971.

\bibitem{DH} David Harel, On visual formalisms. In: J. Glasgow et al
	eds., {\em Diagrammatic Reasoning: Cognitive and Computational
	Perspective} (MIT Press, 1995) 235--271.

\bibitem{CH} Claudio Hermida, {\em From Coherent Structures to Universal
	Properties}, preprint, 1999,
	\texttt{http://www.maths.usyd.edu.au:8000/u/hermida/}.

\bibitem{TL} Tom Leinster, 
	{\em Generalized Enrichment for Categories and Multicategories},
	preprint, 1999, \texttt{math.CT/9901139}.

\bibitem{MM1} Michael Makkai, Avoiding the axiom of choice in general
	category theory, {\em Journal of Pure and Applied Algebra} 108
	(1996) 109--173.

\bibitem{MM2} Michael Makkai, Towards a categorical foundation of
	mathematics. In: J. A. Makowsky and E. V. Ravve eds., {\em Logic
	Colloquium '95}, Lecture Notes in Logic (Springer-Verlag, 1998),
	153--190.

\bibitem{HM} Hiroyuki Miyoshi, Combinatorial structures for
	higher-dimensional categories, talk at CT97.

\bibitem{JP} Jacques Penon, Approche polygraphique des $\infty$-cat\'eogires
	non strictes, {\em Cahiers de Topologie et G\'eom\'etrie
	Diff\'erentielle Cat\'egoriques}, Volume~XL-1 (1999), 31--80.

\bibitem{DP} Detlef Plump, Hypergraph rewriting: critical pairs and
	undicidability of confluence. In: M. R. Sleep, R. Plasmeijer,
	and M. van Eekelen, eds., {\em Term Graph Rewriting, Theroy and
	Practice} (Wiley \& Sons, Chichester, 1993) 201--213.

\end{thebibliography}
\end{document}